\documentclass[%
 aip,
 amsmath,amssymb,
 reprint,%
]{revtex4-1}

\usepackage{graphicx}
\usepackage{dcolumn}
\usepackage{bm}

\usepackage[utf8]{inputenc}
\usepackage[T1]{fontenc}
\usepackage{mathptmx}
\usepackage{etoolbox}
\usepackage{subfigure}
\usepackage{pifont}

\usepackage[colorlinks=true,linkcolor=black,anchorcolor=black,citecolor=black,filecolor=black,menucolor=black,runcolor=black,urlcolor=black]{hyperref}

\usepackage{float}
\makeatletter
\let\newfloat\newfloat@ltx
\makeatother

\usepackage{algorithm}
\usepackage{algpseudocode}
\usepackage{color}
\usepackage{soul}

\DeclareMathOperator*{\argmin}{arg\,min}
\DeclareMathOperator{\grad}{grad}
\DeclareMathOperator{\Hess}{Hess}

\newtheorem{example}{Example}

\makeatletter
\def\@email#1#2{%
 \endgroup
 \patchcmd{\titleblock@produce}
  {\frontmatter@RRAPformat}
  {\frontmatter@RRAPformat{\produce@RRAP{*#1\href{mailto:#2}{#2}}}\frontmatter@RRAPformat}
  {}{}
}%
\makeatother
\begin{document}

\preprint{AIP/123-QED}

\title[Locating saddle points using gradient extremals on manifolds adaptively revealed as point clouds]{Locating saddle points using gradient extremals on manifolds adaptively revealed as point clouds}
\author{A. Georgiou}
\email{ageorgi3@jhu.edu}
\affiliation{ 
    Department of Chemical and Biomolecular Engineering, Johns Hopkins University, Baltimore, Maryland 21218, USA}%

\author{H. Vandecasteele}%
\affiliation{ 
Department of Computer Science, KU Leuven, 3001 Leuven, Belgium
}%

\author{J. M. Bello-Rivas}%
\affiliation{ 
Department of Applied Mathematics \& Statistics, Johns Hopkins University, Baltimore, Maryland 21218, USA
}%

\author{I. Kevrekidis}
\email{yannisk@jhu.edu}
\affiliation{ 
    Department of Chemical and Biomolecular Engineering, Johns Hopkins University, Baltimore, Maryland 21218, USA}%
\affiliation{ 
Department of Applied Mathematics \& Statistics, Johns Hopkins University, Baltimore, Maryland 21218, USA
}%

\date{\today}

\begin{abstract}
Steady states are invaluable in the study of dynamical systems. High-dimensional dynamical systems, due to a separation of time-scales, often evolve towards a lower dimensional manifold $M$. We introduce an approach to locate saddle points (and other fixed points) that utilizes gradient extremals on such \textit{a priori} unknown (Riemannian) manifolds, defined by adaptively sampled point clouds, with local coordinates discovered on-the-fly through manifold learning. The technique, which efficiently biases the dynamical system along a curve (as opposed to exhaustively exploring the state space), requires knowledge of a single minimum and the ability to sample around an arbitrary point. We demonstrate the effectiveness of the technique on the M\"{u}ller-Brown potential mapped onto an unknown surface (namely, a sphere). Previous work employed a similar algorithmic framework to find saddle points using Newton trajectories and gentlest ascent dynamics; we therefore also offer a brief comparison with these methods.
\end{abstract}

\maketitle

\begin{quotation}
The energy landscapes for molecular dynamics, protein folding, glassy materials, etc. may exhibit large numbers of minima and transition states. Locating these critical points---especially in high dimensions---is challenging. Current approaches (e.g. nudged elastic band, string method) start with two minima and \textit{a priori} knowledge of collective variables which represent the system well in a lower dimensional setting; other methods attempt to exhaustively explore the state space (metadynamics, adaptive biasing force, iMapD). In this paper, we present a technique to systematically bias dynamical systems from a single minimum to a saddle point (transition state) without \textit{a priori} knowledge of good collective variables.
\end{quotation}

\section{\label{sec:intro}Introduction}

Locating saddle points of dynamical systems is of primary importance to a variety of applications, including the location of transition states of chemical systems described at the atomic level~\cite{Bochevarov2013} or the study of transitions between molecular configurations in molecular dynamics.~\cite{leimkuhler2015}
Unfortunately, finding these transition states simply via simulation is time-prohibitive: the trajectories spend the majority of their time in a basin of attraction around a stable equilibrium, rarely traveling to a different basin of attraction (potential well). Often, trajectories of the high-dimensional system quickly evolve to, and then move slowly along, a low-dimensional manifold (for example, due to a time-scale separation). The coarse (collective) variables that parameterize the slow manifold are either known \textit{a priori} or can be determined through manifold learning techniques.

Many techniques exist to locate saddle points of these systems. The metadynamics~\cite{laio2002} and adaptive biasing force~\cite{darve2008} methods attempt to exhaustively explore the whole space, yet they require \textit{a priori} knowledge of collective variables; iMapD~\cite{chiavazzo2017} alleviates this need by exploiting manifold learning techniques. Another class of methods, namely nudged elastic band or string method variants,~\cite{olender1996,Jonsson1998,Weinan2002} trace a curve between two known minima: specifically, they minimize the energy at each point along a path connecting two minima (reactant and product) while attempting to maintain equidistant spacing between the points. These methods require a global set of collective variables. A third class of methods require only a single minimum, utilizing local information to follow a continuous curve to a nearby fixed point (hopefully a saddle). These one-dimensional curves include Newton trajectories,~\cite{quapp1998,hirsch2004} gentlest ascent dynamics~\cite{Weinan2010} or dimer method~\cite{henkelman1999} trajectories, and gradient extremals,~\cite{basilevsky1982} which we study here. These methods still require a global set of collective variables and were originally formulated for Euclidean space.

Gradient extremals have been extended to smooth Riemannian manifolds to locate fixed points.~\cite{Rowe1982,Filippidis2013} In this work, we forego \textit{a priori} knowledge of good collective variables as well as an explicit formulation of the manifold. Instead, we follow gradient extremals on initially unexplored Riemannian manifolds gradually revealed by adaptively sampled point clouds. We utilize local sets of collective variables that are discovered on-the-fly through manifold learning (here, diffusion maps~\cite{COIFMAN2006}). By sampling along the manifold, discovering a good local chart, resolving the path locally on the chart, and repeating this process at the newly found location, the system is driven from an initial conformation (typically a stable equilibrium) to a saddle point. Our methodology successfully locates new fixed points using a single initial point and without need for \textit{a priori} knowledge of collective variables. The method does not exhaustively explore the state space, but rather follows a single path at a time. In prior work we utilized a similar algorithmic approach, but with the resolved paths being those of Newton trajectories~\cite{Bello-Rivas2023} and gentlest ascent dynamics.~\cite{bellorivas2023gentlest}

\section{\label{sec:GE}Gradient Extremals}

Gradient extremals are a type of curve that connects fixed points of a dynamical system. We will assume, because of our motivating examples, that the dynamical system is a gradient system, i.e., the vector field $X(x) = -DE(x)$ for some twice-differentiable potential $E$. This assumption can be relaxed to arbitrary vector fields by working with the squared magnitude of the vector field.~\cite{lucia2002} An example of this can be found towards the end of Section \ref{sec:grad-ex-euc}. Gradient extremals can be thought of as following the ``valleys, ridges, cirques, or cliffs'' of the potential $E$.~\cite{hoffman1986}

\subsection{\label{sec:grad-ex-euc}Gradient Extremals in Euclidean Space}

\begin{figure*}
\centering
\begin{subfigure}[]
  \centering
  \includegraphics[trim={3cm 7cm 2cm 6cm},clip,width=0.55\linewidth]{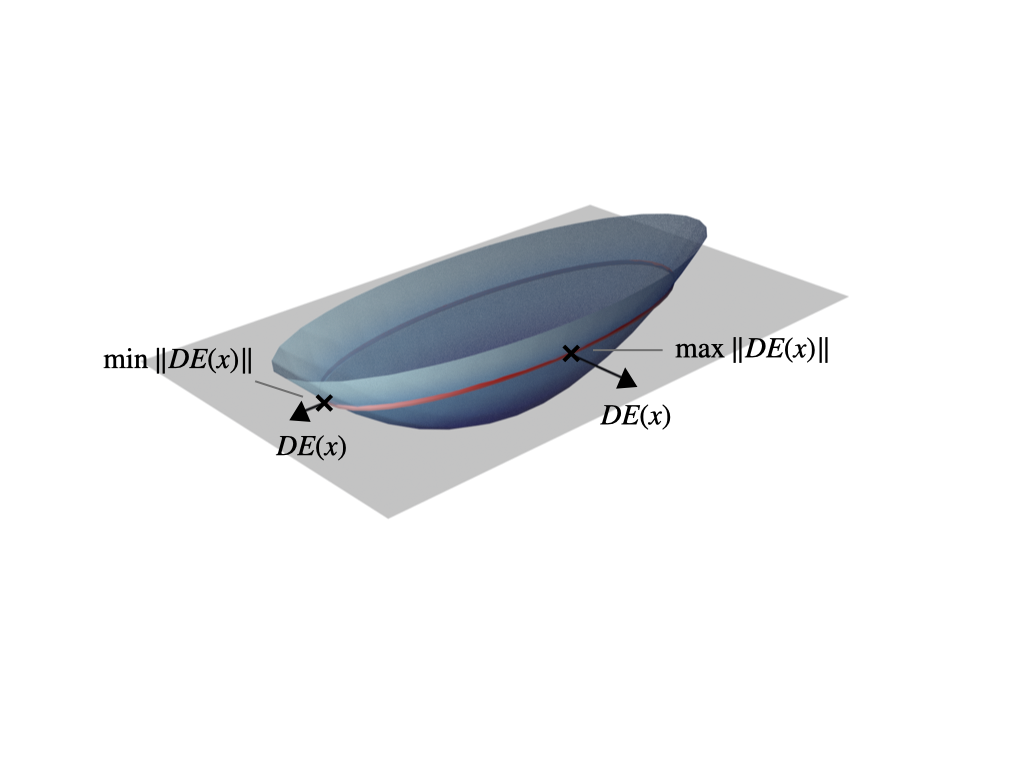}
  \label{fig:subd}
\end{subfigure}%
\begin{subfigure}[]
  \centering
  \includegraphics[trim={0cm 0cm 0cm 0cm},clip,width=0.35\linewidth]{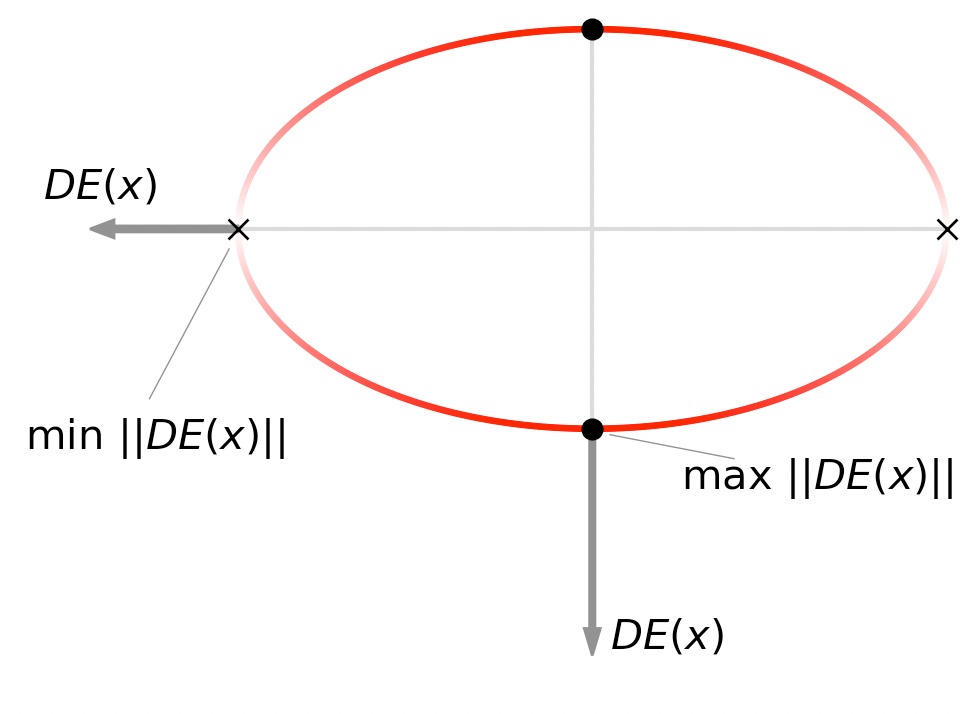}
\end{subfigure}
\caption{\textbf{(a)} A harmonic potential $E$ intersected by a plane at $E=L$. The elliptical contour line is colored by the length of the corresponding force, with red being a maximum and white a minimum. The points of the gradient extremal corresponding to the min and max of $\|DE\|$ at $E=L$ are marked with \ding{54}. \textbf{(b)} A 2D schematic of (a).}
\label{fig:schematic}
\end{figure*}

\begin{figure}
\centering
    \includegraphics[trim={0cm 0cm 0cm 0cm},clip,width=1\linewidth]{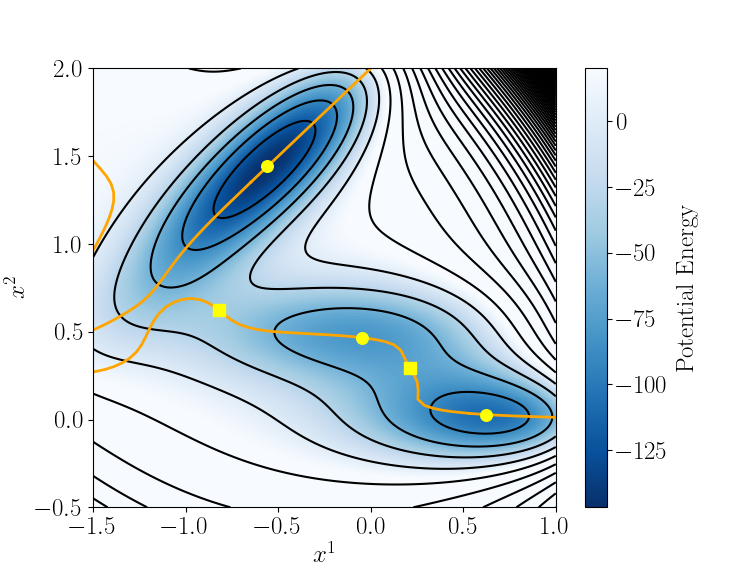}
\caption{Gradient extremals (following the eigenvector corresponding to the smallest eigenvalue of the Hessian of the potential) on the M\"{u}ller-Brown potential, which is defined in~\eqref{eq:muller-brown}. Mimima are depicted by yellow circles; saddle points by squares. Note that the top minimum is not directly connected to the bottom two. They are instead connected through a long detour~\cite{ohno2004} that extends beyond the boundaries of the figure.}
\label{fig:GE_MB}
\end{figure}

\begin{figure}
\centering
\begin{subfigure}[]
  \centering
  \includegraphics[width=0.9\linewidth]{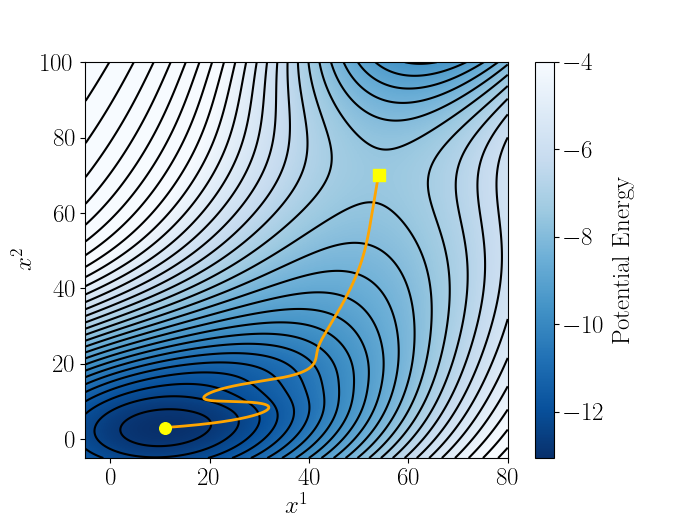}
  \label{fig:sube}
\end{subfigure}%
\begin{subfigure}[]
  \centering
  \includegraphics[trim={0.7cm 0cm 2cm 2cm},clip,width=0.9\linewidth]{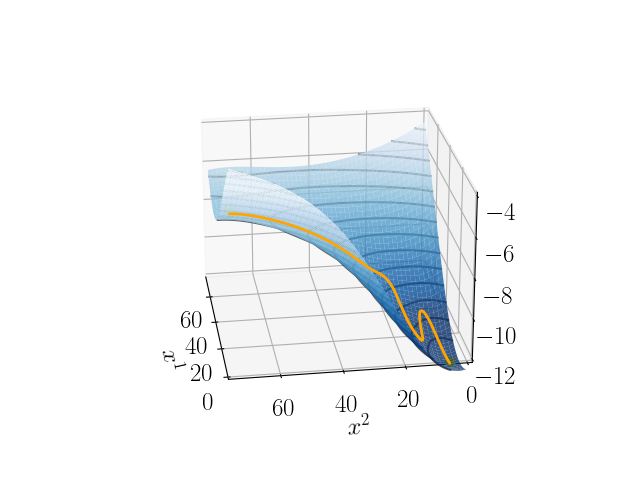}
\end{subfigure}
\caption{The gradient extremal following the eigenvector corresponding to the smallest eigenvalue of the Hessian of~\eqref{eq:yannik}. Note how the gradient extremal ``meanders,'' or exhibits turning points.}
\label{fig:yannik}
\end{figure}

Gradient extremal curves~\cite{basilevsky1982, hoffman1986, quapp1989, lucia2002} are defined as the locus of points that minimize (or maximize) the gradient norm along each level set of the potential $E$
\begin{equation}
  \label{eq:lucia_optimization}
  \underset{x\in \mathbb{R}^n}{\argmin} (\max) \{ \| D E(x) \|^2 : E(x)=L\}.
\end{equation}
The points on the level set highlighted in Figure~\ref{fig:schematic} 
that minimize (respectively, maximize) the gradient are the farthest (respectively, closest) from the fixed point. Equivalently, one may formulate gradient extremals as the locus of points where the gradient of the potential is an eigenvector of the Hessian $D^2 E(x)$,
\begin{equation} \label{eq:classical_ge}
    D^2 E(x) DE(x) = \lambda DE(x),
\end{equation}
along each level set $E(x)=L$. Here $\lambda$ is the associated eigenvalue.
One can verify that both methods are equivalent by applying the method of Lagrange multipliers on~\eqref{eq:lucia_optimization}.
Gradient extremals can be followed through continuation methods~\cite{continuation2003} where the level set value, $L$, acts as the parameter of interest and~\eqref{eq:lucia_optimization} is optimized for each $L$.~\cite{lucia2004} Note that a specific gradient extremal curve may not traverse all fixed points---it might instead connect just a subset of them.~\cite{hirschquapp2004} See Figure~\ref{fig:GE_MB} 
for an example of gradient extremals on the M\"{u}ller-Brown potential (defined in ~\eqref{eq:muller-brown}).
We also note that a gradient extremal can exhibit turning points (or  ``meander''), as seen in Figure~\ref{fig:yannik},
which depicts a gradient extremal on the potential defined by
\begin{align}
    \label{eq:yannik}
E&(x^1, x^2)= \nonumber \\
&-\sum_{i = 1}^4 A_i \exp\left(a_i (x_i^1-x^1_{0i})^2 + b_i(x^2_i-x^2_{0i})^2 + c_i(x^1_i x^2_i)\right)
\end{align}
with the coefficients listed in Table~\ref{tab:Yannik}. The superscripts here, and throughout the paper, \textit{refer to the component indices} (Einstein notation) unless otherwise noted.

\begin{table}
\caption{\label{tab:Yannik}Coefficients of the potential defined by \eqref{eq:yannik}.}
\begin{ruledtabular}
\renewcommand{\arraystretch}{1.25}
\begin{tabular}{c c c c c c c}
 $i$ & $A_i$ & $a_i$ & $b_i$ & $c_i$ & $x^1_{0i}$ & $x^2_{0i}$\\ [0.5ex]
 \hline
 1 & 10 & $\frac{-1}{3000}$ & $\frac{-1}{1000}$ & 0 & 0 & 0\\
 2 & $-0.4$ & $\frac{-1}{300}$ & $\frac{-1}{30}$ & $\frac{1}{200}$ & 15 & 10 \\
 3 & 0.8 & $\frac{-1}{1500}$ & $\frac{-1}{1500}$ & $\frac{1}{2000}$ & 25 & 100 \\
 4 & 6 & $\frac{-1}{3000}$ & $\frac{-1}{3000}$ & $\frac{3}{20000}$ & 40 & 30\\ [1ex]
\end{tabular}
\end{ruledtabular}
\end{table}

\begin{figure}
\centering
  \includegraphics[trim={0cm 0cm 0cm 0cm},clip,width=0.9\linewidth]{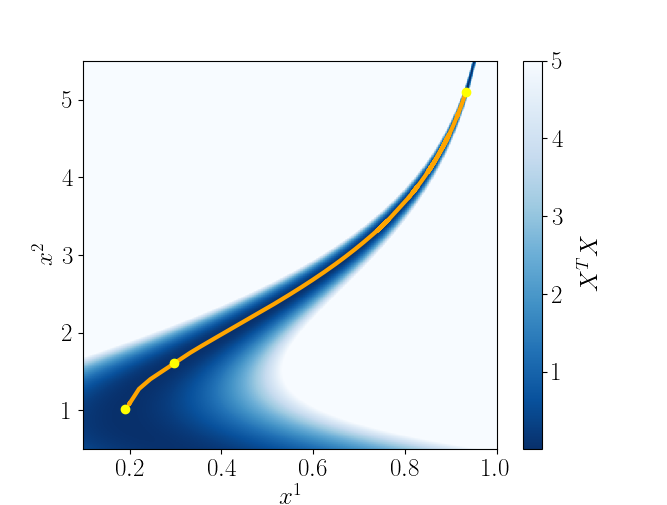}
\caption{\label{fig:CSTR}Gradient extremal (orange curve) on the graph of the squared magnitude, $X^TX$, for the CSTR. The steady states are shown in yellow.}
\end{figure}

We previously assumed that the dynamical system is of gradient type. However, this assumption can be removed by replacing the potential by the squared magnitude of any given vector field $X$.~\cite{lucia2002} The gradient extremal then becomes the locus of points that solve the constrained optimization problem
\begin{equation} \label{eq:lucia_optimization_general}
\underset{x \in \mathbb{R}^n}{\argmin}(\max) \ \{ \lVert D (X^TX) \rVert^2: X^TX=L\}.
\end{equation}
Notice that $E$ in~\eqref{eq:lucia_optimization} is replaced by the squared magnitude $X^TX$. We illustrate this technique with the vector field describing a first order, exothermic, irreversible reaction $A\rightarrow B$ taking place in a continuous-stirred tank reactor (CSTR). The governing mass and energy balances are
\begin{align}
    \label{eq:CSTR}
\frac{dx^1}{dt} &= -x^1 + D_a(1-x^1)\exp{\big(x^2\big)} \nonumber \\ 
\frac{dx^2}{dt} &= -x^2+BD_a(1-x^1)\exp{\big(x^2\big)}-\beta(x^2-x^2_{c})
\end{align}
where $x^1$ is the conversion of species $A$ and $x^2$ is dimensionless temperature as detailed in Ref.~\onlinecite{Uppal1974}.  (Please remember that the superscripts in $x^1$ and $x^2$ are indices, not powers, in Einstein notation.) With parameters set to $D_a = 0.085$, $B=22$, $\beta=3$, $x^2_{c}=-0.04$, the dynamical system exhibits two unstable steady states, one stable steady state, and one stable limit cycle. Figure~\ref{fig:CSTR} 
illustrates gradient extremals on the squared magnitude of the vector field~\eqref{eq:CSTR}.

\subsection{\label{sec:GE_Riemmanian}Gradient Extremals on Riemannian Manifolds}
We now turn our attention to resolving gradient extremals for vector fields defined on a manifold. Recall that the purpose of the proposed algorithm is to connect fixed points (hopefully a local minimum to a kinetically relevant nearby saddle): if one could sample and map a chart (defined as a small portion of the manifold represented with reduced coordinates) that includes both minimum and saddle point, one could directly solve for a gradient extremal by treating the chart as a Euclidean space. The resulting path would not be a true gradient extremal, but would still properly link the two critical points, which was our original intent. However, chances are that multiple local charts will be needed, whether because (a) the curvature of the manifold makes the projection onto a single chart impossible or (b) a single chart does not contain both the initial and final fixed points. In the former case, increased curvature makes an accurate, invertible mapping between manifold and local chart difficult: choosing multiple smaller areas to work with alleviates this issue. In the latter case, unless both connected fixed points lie on the chart, a gradient extremal would terminate at the boundary of the chart without connecting to the next fixed point of interest.

To ensure the path transitions properly from one chart to the next, we use differential geometry to resolve the gradient extremal in the Riemannian setting of a $d$-dimensional smooth manifold, $M \subset \mathbb{R}^n$, with Riemannian metric $g$. 
The Riemannian metric $g$ defines the notions of distances and angles on the manifold by assigning inner products on the tangent spaces of the manifold. Let $E : M \to \mathbb{R}$ be the potential energy function on the manifold. We create a mapping $\varphi: U \subset M \to \mathbb{R}^d$ that defines a system of coordinates of an open subset $U \subset M $ on a local chart in $\mathbb{R}^d$. Let $\psi : \mathbb{R}^d \to U$, $u \mapsto \psi(u)$, be the inverse map of the chart $\varphi$. We define a Hessian matrix in terms of the covariant derivative on $(M, g)$ with respect to the Levi-Civita connection. 
A connection describes how a vector changes as it moves from one point on a manifold to the next.
We use the following notation when working on $(M,g)$, with $Z=\psi^{\star} E: \mathbb{R}^d \to \mathbb{R}$ now representing the pullback energy:
\begin{equation*}
     Z=\psi^{\star} E = E\circ \psi. 
\end{equation*}
Thus, the Euclidean formulation of gradient extremals \eqref{eq:classical_ge} now becomes, on a local chart,
\begin{equation}\label{eq:hessian-eigenproblem}
  \Hess Z \, \grad Z = \lambda \grad Z.
\end{equation}
First, we define our Riemannian metric, $g$, as
\begin{equation}\label{eq:metric}
  g = \sum_{i, j = 1}^d g_{ij} \, \text{d}\varphi^i \otimes \text{d} \varphi^j,
\end{equation}
where $g_{ij} = (D \psi^\top D \psi)_{ij}$ for $i, j = 1, \dotsc, d$. This metric allows us to define an inner product between two tangent vectors $S$ and $T$ defined in local coordinates: $g(S,T) = \sum^d_{i,j=1} g_{ij}S^iT^i$. We can now compute the gradient and Hessian on the manifold. Each component of the gradient becomes
\begin{equation*}
  \grad(Z)^i = \sum^d_{j=1} g^{ij}\frac{\partial Z}{\partial \varphi^j},
\end{equation*}
where $g^{ij}$ refers to element $(i,j)$ of the inverse metric tensor. The Hessian is defined in terms of the covariant derivative:
\begin{equation*}
  \Hess(Z)T
  =
   \nabla_{T} \grad({Z})
\end{equation*}
where the covariant derivative, which allows us to specify the parallel transport along tangent vectors of a manifold, is described as
\begin{equation*}\label{eq:riemannian-parallel-transport}
  \nabla_{T} Y
  =
  \sum_{k = 1}^d \left\{
    \sum_{i = 1}^d
    \left(
      \frac{\partial Y^k}{\partial \varphi^i}
      +
      \sum_{j = 1}^d \Gamma_{ij}^k Y^j
    \right)
    {T^i}
  \right\}
  \frac{\partial}{\partial\varphi^k}.
\end{equation*}
The above expression constitutes the covariant derivative of $Y$ in the direction of $T$. The coefficients $\Gamma_{ij}^k$ act as a correction to otherwise straight trajectories on the tangent plane so that they remain on the manifold. These coefficients are unique to a given affine connection. Here we use the Levi-Civita connection, the unique symmetric connection compatible with the metric, for which the coefficients, known as Christoffel symbols, are given by
\begin{equation*}
  \Gamma^\ell_{jk}
  =
  \frac{1}{2}\sum_{i = 1}^d g^{\ell i} \left( \frac{\partial  g_{ij}}{\partial \varphi^k} + \frac{\partial  g_{ik}}{\partial \varphi^j} - \frac{\partial g_{jk}}{\partial \varphi^i} \right).
\end{equation*}
Beginning with an initial point $p_0$, its energy level $L = E(p_0)$, and the smallest (or largest) eigenvalue, $\lambda$, of the Hessian of the pullback energy evaluated at the initial point, we can now use continuation to solve the following system of $d+1$ equations:
\begin{equation*}
  \left\{
    \begin{aligned}
      &Z - L = 0 \\
      &\Hess(Z) \grad(Z) - \lambda\ \grad(Z) = 0
    \end{aligned}
  \right.
  \label{eq:continuation}
\end{equation*}
Note that naively using the Euclidean formulation of gradient extremals as defined in~\eqref{eq:classical_ge} on the local chart would be erroneous because it fails to take into account the geometry of the manifold. Indeed, we would obtain curves that would not, in general, coincide with the true gradient extremals. This is akin to how a straight line joining New York and Madrid drawn on a (Mercator) map is not actually a ``straight-line'' on the globe (\emph{i.e.,} a segment of a great circle ---also known as a geodesic).

\begin{figure}
\centering
\begin{subfigure}[]
  \centering
    \includegraphics[trim={18cm 10cm 18cm 10cm},clip,width=0.6\linewidth]{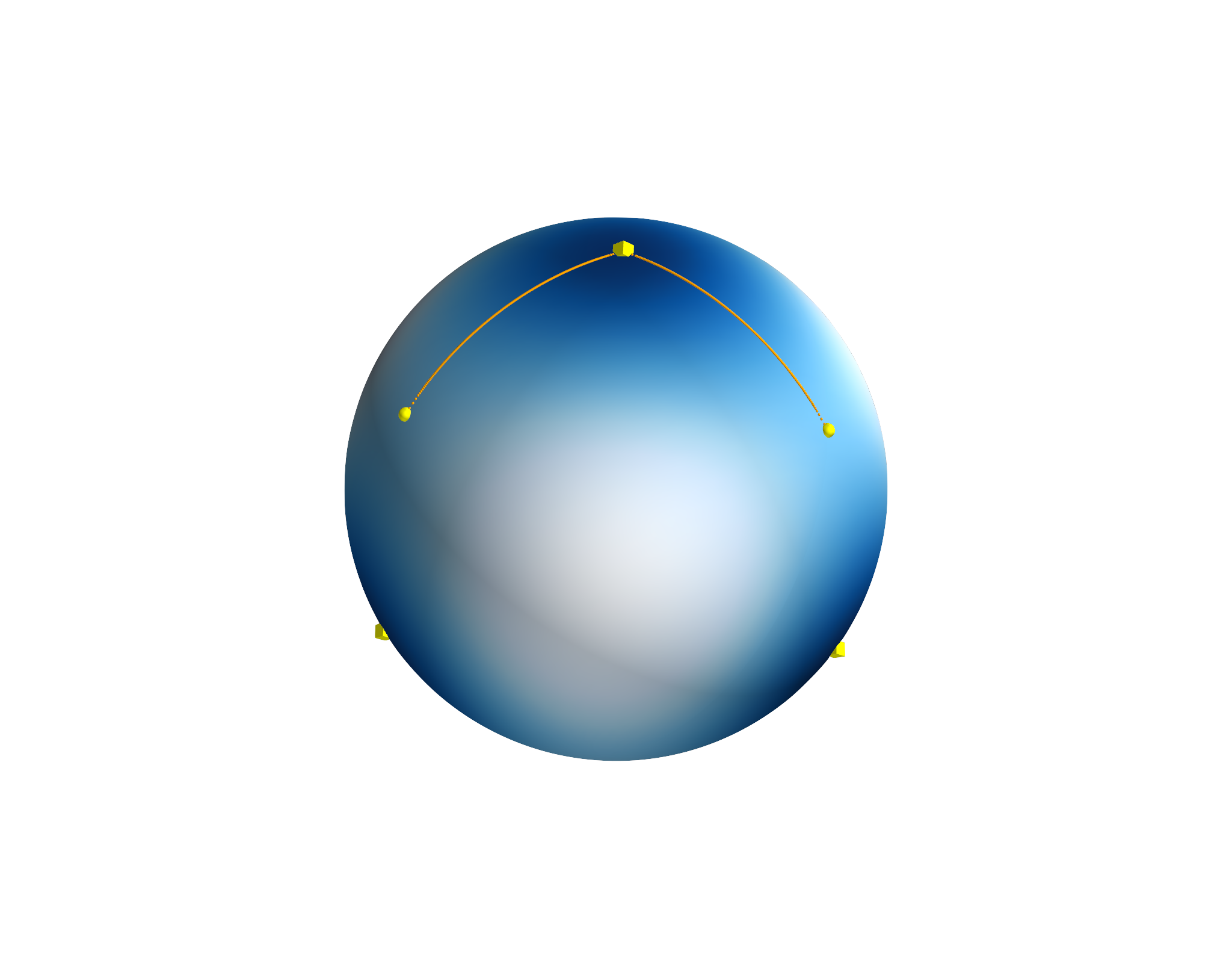}
  \label{fig:subf}
\end{subfigure}%
\begin{subfigure}[]
  \centering
\includegraphics[width=0.9\linewidth]{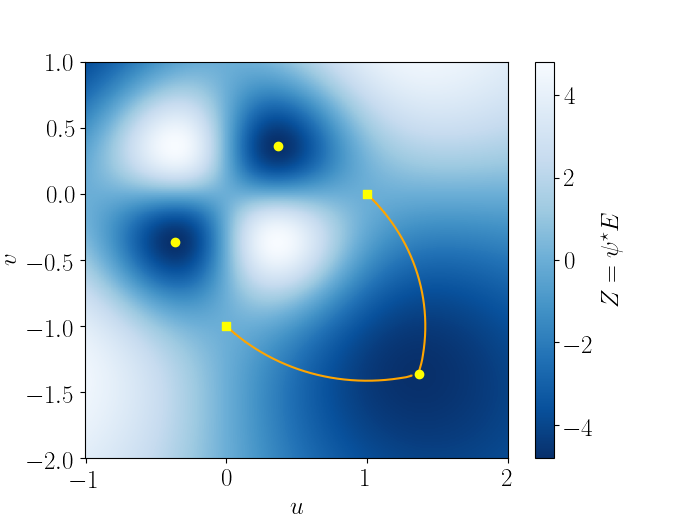}
\end{subfigure}
\caption{\textbf{(a)} Two gradient extremal curves shown on the sphere $\mathbb{S}^2 = \{ (x,y,z) \in \mathbb{R}^3 \mid x^2+y^2+z^2=1\}$ with the potential $E(x,y,z)=xyz$. \textbf{(b)} The stereographic projection from the North Pole to the tangent plane at the South pole is shown along with the associated pullback potential energy $\psi^\star E$. Full calculation detailed in Example~\ref{ex:zorro}.}
\label{fig:Zorro}
\end{figure}

\begin{example}
We illustrate a gradient extremal resolved on a smooth Riemannian manifold defined by the unit sphere $\mathbb{S}^2 = \{ (x,y,z) \in \mathbb{R}^3 \mid x^2+y^2+z^2=1\}$ with potential energy $E(x,y,z)=xyz$ in Figure~\ref{fig:Zorro}. 
We define the mapping $\varphi$ as the stereographic projection from the North pole to the tangent plane at the South pole such that $\varphi(x,y,z) = \left(\frac{x}{1-z}, \frac{y}{1-z}\right) \equiv (u,v)$ and its inverse as $\psi(u, v) = \frac{1}{u^2+v^2+1}\left(2u,2v,u^2+v^2-1\right)$, and the pullback energy as $Z(u,v) = \frac{4uv\left(u^2+v^2-1\right)}{\left(u^2+v^2+1\right)^3}$. All fixed points except the North pole are mapped onto a single chart; in our actual numerical scheme, we relax this assumption to allow for sampling of only a small portion of the manifold at a time. We refer the reader to Ref.~\onlinecite{bellorivas2023gentlest} for a full derivation of the gradient and Hessian used in this example.
\label{ex:zorro}
\end{example}

As we did in Section \ref{sec:grad-ex-euc}, we can relax the assumption that $E$ corresponds to the potential energy of a gradient vector field if we work instead with the squared magnitude of a general vector field $X$ on the manifold. If we let $\ell(u)$ be the pullback of the squared magnitude evaluated at point $u$, then the counterpart to~\eqref{eq:lucia_optimization_general} on a Riemannian manifold $(M, g)$ becomes
\begin{equation*}
  \argmin_{u \in \ell^{-1}(\{ L \})}(\max) \left\{ g_u(\grad \ell(u), \grad \ell(u)) \right\},
\end{equation*}
where
\begin{align*}
  \ell \colon M &\to \mathbb{R}_{\ge 0} \\
  u & \mapsto \ell(u) = g_u(\psi^{\star}X, \psi^{\star}X) = \psi^{\star}(X^TX).
\end{align*}
Here, $\grad$ is again the Riemannian gradient and  $\ell(u)$ corresponds to the pullback of the term $X^\top X$ in the Euclidean case.

\begin{figure}
\centering
\begin{subfigure}[]
  \centering
  \includegraphics[width=0.85\linewidth]{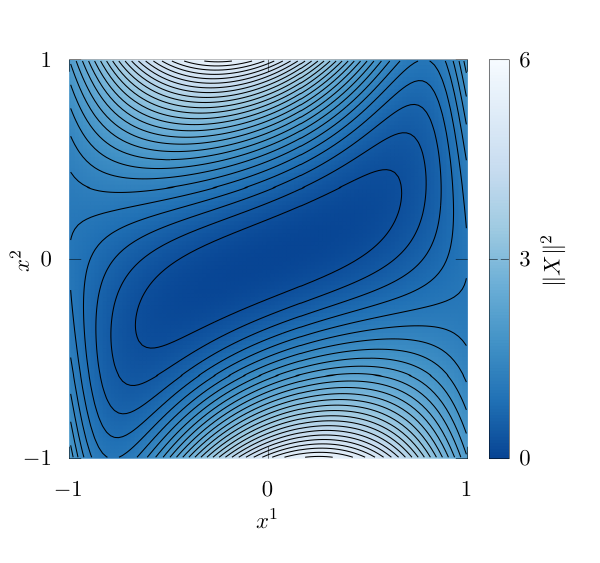}
  \label{fig:suba}
\end{subfigure}%
\begin{subfigure}[]
  \centering
  \includegraphics[trim={0.0cm 0cm 0cm 0cm},clip,width=0.9\linewidth]{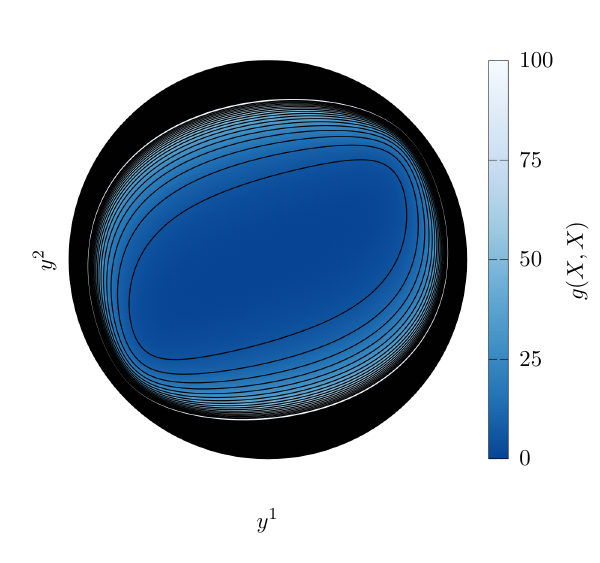}
\end{subfigure}
\caption{Squared length function of the van der Pol vector field ($\mu = 2$) on two manifolds. \textbf{(a)} Euclidean plane. \textbf{(b)} Poincar\'e disk model of the hyperbolic plane.}
\label{fig:squared-length-function}
\end{figure}

As an example of a non-gradient system, let us study the van der Pol oscillator on a Riemannian manifold. First, let $\mathbb{R}^2$ be the Euclidean plane with coordinates $(x^1, x^2)$ and the usual Euclidean metric.
The van der Pol oscillator is given by the vector field $X$:
\begin{align}
    \frac{dx^1}{dt} &= x^2  \nonumber \\
    \frac{dx^2}{dt} &= \mu(1-(x^1)^2)x^2 -x^1 \nonumber
\end{align}
for some constant $\mu \in \mathbb{R}$.
Next, the Poincar\'e disk model of the hyperbolic plane is the Riemannian manifold given by
\begin{equation*}
  \mathbb{H}^2
  =
  \left\{
    (y^1, y^2) \in \mathbb{R}^2 \, : \,
    (y^1)^2 + (y^2)^2 < 1
  \right\}
\end{equation*}
and the metric
\begin{equation*}
  g
  =
  \frac{4 (\mathrm{d}y^1 \otimes \mathrm{d} y^1 + \mathrm{d}y^2 \otimes \mathrm{d} y^2)}{\left( 1 - (y^1)^2 - (y^2)^2 \right)}.
\end{equation*}
We define the van der Pol vector field $X$ on $\mathbb{H}^2$ simply by restricting the Euclidean vector field to the subset $\mathbb{H}^2$.
In other words, we retain only the vector field $X$ as defined in the disk $\mathbb{H}^2$ and this yields, in the new hyperbolic metric $g$ above, \textit{a different squared length function} than the Euclidean one.
Figure~\ref{fig:squared-length-function} shows the difference between the squared length function of the van der Pol oscillator $X$ in $\mathbb{R}^2$ and its restriction to $\mathbb{H}^2$.

\section{\label{sec:Algorithm}Algorithm}
We consider a high-dimensional dynamical system $X(x) \in \mathbb{R}^n$ whose trajectories become quickly attracted to a slow $d$-dimensional manifold $M \subset \mathbb{R}^n \ (n \gg d)$. Suppose we have an initial point $p_0 \in M$ near a minimum, obtained by running a single long-enough gradient descent trajectory, such that the fast dynamics have subsided, and movement on the manifold dominates. Next, we draw $N$ samples in a neighborhood $U \subset M$ around $p_0$. Any sampling technique appropriate to the application can be used; in the subsequent example, we randomly perturb $p_0$ in ambient space and run the $N$ trajectories for a short time $\tau$, such that they now effectively lie on the manifold. Another option would be to use umbrella sampling~\cite{Torrie1977,fiorin2013} based on the reduced local coordinates.\cite{Bello-Rivas2023}

Next, we apply a dimensionality reduction technique such as diffusion maps~\cite{COIFMAN2006} to uncover a set of reduced, latent coordinates that describe the manifold in the neighborhood of the initial point. It is important to emphasize that we use \textit{global} manifold learning techniques on a neighborhood on the manifold: as a result, we obtain a good set of coordinates describing the single, \textit{local} chart. A different set of collective variables may be considered for each subsequent chart, which could be better suited for the associated local neighborhood (which could correspond to a different operating regime, stage, conformation, etc. of the system). We then use a Gaussian process to learn the local mappings between the ambient space and discovered chart $\varphi: U \rightarrow \mathbb{R}^d$ and its inverse $\psi: \mathbb{R}^d \rightarrow U $, as well as a mapping of the reduced coordinates to the potential energy, the pullback energy  $Z: \mathbb{R}^d \rightarrow \mathbb{R}$. In the simple examples of this paper, the potential function in ambient space is given. In more complex applications (e.g. effective descriptions of molecular dynamics simulations in terms of collective variables), where this potential function is not explicitly known, some additional local computations (e.g. through umbrella sampling) would be required to estimate its gradient (mean force) at each step. Further, in the case of a high-dimensional ambient space (of higher dimensions than the examples shown here), a Gaussian process may inadequately lift from latent to ambient space. Rather than fitting a function for $\psi$ in these more complex applications, Ref.~\onlinecite{Bello-Rivas2023} presents instead a method using biased sampling by simulating a stochastic differential equation to estimate $\psi(u)$.

After learning the mappings between the manifold and the local chart, we can now compute the gradient extremal in the Riemannian setting as explained in Section~\ref{sec:GE_Riemmanian}. Here, we use a pseudo-arclength continuation method~\cite{continuation2003} with adaptive step sizing, stopping when an equilibrium point or chart boundary is reached. Detection of a fixed point can be done by monitoring the norm of the vector field (either using the vector field itself in ambient space, or by inspecting the gradient of the pullback energy on the chart). If using a Gaussian Process to obtain the diffeomorphisms, detection of a chart boundary can be accomplished by monitoring that the norm of the covariance matrix at a given point remains below a chosen threshold. If a boundary is reached, it is necessary to switch charts --- that is, to resample around the point at the boundary of the previous chart and relearn the mappings and pullback energy function. Care has to be taken to continue the curve in the same direction between local charts, potentially oriented in opposing directions; this can be accomplished by  mapping the last few ambient points of the curve onto the new local chart and ensuring continuation occurs in the same direction (e.g. the angle between the previous curve mapped onto the new chart and the continuation direction is acute). The full methodology is summarized in Algorithm~\ref{alg:ge}.

\begin{algorithm*}
\renewcommand{\thealgorithm}{ 1}
  \caption{Continuation of gradient extremal curves on manifolds.}\label{alg:ge}
  \begin{algorithmic}
    \Require initial point $p=p_0$ near minimum, samples per iteration $N$, continuation parameters, threshold $\rho > 0$ used in the convergence criterion.
    \Ensure equilibrium $p_\star$ of $X$.
    \For{$n = 1, 2, 3, \dotsc$}
      \State Sample $N$ points from a neighborhood $U \subset M$ of $p_{n\text{-}1}$.
      \State Use manifold learning to obtain local coordinates $\varphi \colon U \rightarrow \mathbb{R}^d$
      \State Obtain a parameterization $\psi \colon V \rightarrow U$, where $V = \varphi(U) \subseteq \mathbb{R}^d$.
      \State Approximate $Z = \psi^\star E$ via Gaussian process regression.
      \State{$q \leftarrow \varphi(p_{n-1})$}
      \If{n = 1}
        \State{$v \leftarrow$ Random direction in $\mathbb{R}^d$}
      \Else
        \State{$v \leftarrow \varphi_\star(q) w_{n-1}$}
      \EndIf
      \State $q_n, v_n \leftarrow \textsc{resolve\_extremal\_curve}(q, v)$ \Comment{See Algorithm~\ref{alg:ge-aux}.}
      \State $p_n, w_n \leftarrow \psi(q_n), \psi_\star(q_n) v_n$
      \If{$\ \| \grad Z(q_n) \| < \rho$ and $n>1$}
        \State \Return $p_\star \leftarrow p_n$
      \EndIf
    \EndFor
  \end{algorithmic}
\end{algorithm*}

\begin{algorithm*}
\renewcommand{\thealgorithm}{ 2}
  \caption{\textsc{resolve\_extremal\_curve}}\label{alg:ge-aux}
  \begin{algorithmic}
    \Require Initial point $q_0 \in \mathbb{R}^d$ and velocity $v_0 \in \mathbb{R}^d$.
    \Ensure Final point $q$ and velocity $v$ of the extremal curve $\gamma$.
    \State Construct numerical continuation curve $\gamma$ for 
    \begin{equation*}
    \left\{ 
    \begin{aligned}
    & \Hess Z\grad Z = \lambda \grad Z \\
    & Z(\gamma(t)) = L
    \end{aligned}
    \right.
    \end{equation*}
    such that $\gamma(0) = q_0$ and $\dot{\gamma}(0) = v_0$. $L$ is the value of the level set at the current time step.
    \State Prolong $\gamma$ until the first instant $t > 0$ such that $\gamma(t) \in V$ but $\gamma(t + s) \not\in V$ for any $s > 0$.
    \State \Return $\gamma(t)$, $\dot{\gamma}(t)$
\end{algorithmic}
\end{algorithm*}

While we believe our algorithm is as parsimonious as possible (i.e. starting from a single minimum, using a set of reduced collective variables discovered on-the-fly, and efficiently exploring in a single direction), it is not without limitations. Here, we detail those nuances and limitations:
\begin{description}
\item[Exploration of all directions] Given a particular minimum, there exist multiple different directions along which one can initialize a gradient extremal curve. For example, from a minimum on a 2D manifold, one could follow both sides (positive and negative) of both eigenvectors of the Hessian. Algorithm~\ref{alg:ge} assumes that we always follow a single, pre-chosen direction. A bookkeeping strategy would be required to systematically explore all gradient extremals of a dynamical system. The Global Terrain Algorithm~\cite{lucia2002} and the Navigation Algorithm~\cite{Filippidis2013} utilize a graph structure to conduct such bookkeeping. These techniques could be integrated into Algorithm~\ref{alg:ge} to provide an orderly exploration of the fixed points of the dynamical system.
\item[Disconnected minima] Unfortunately, a single gradient extremal curve is not guaranteed to connect all fixed points of a potential energy surface. Two minima may occasionally lie on separate gradient extremal curves~\cite{hirschquapp2004} or may be connected by a very long detour.~\cite{ohno2004} This is the case in the M\"{u}ller-Brown potential, where the topmost minimum follows a long detour to the bottom two (see Figure~\ref{fig:GE_MB}). An alternative to following the long detour would be to perform a steepest descent trajectory from both sides of the newly located leftmost saddle. Discovery of all saddle points given a single starting minimum is not guaranteed.
\item[Continuation parameters] Our implementation hinges on a predictor-corrector scheme for numerical continuation. Care must be taken to set reasonable tolerance values and step size maxima. This is especially true near fixed points, which give rise to numerical difficulties. Continuation sometimes becomes numerically challenging, e.g. in the neighborhood of a new fixed point in the current chart, because the potential energy surface is locally flat there. To locate the fixed point once we approach its neighborhood, Ref.~\onlinecite{lucia2002} suggests switching to an equation-solving technique such as quadratic acceleration or conjugate gradients to accurately locate the saddle. These techniques can be employed using the learned mappings in the chart of interest.
\item[Sampling and manifold learning] An adequate approximation of the manifold is critical. Poor sampling or failed manifold learning will undermine the success of the algorithm. The more curvature a manifold exhibits, the more challenging it is to create an accurate, invertible mapping between manifold and local chart. Working instead with multiple, smaller charts will alleviate this issue.
\end{description}

\section{Comparison of Potential Curves Connecting Fixed Points}
The briefly aforementioned paths (gradient extremals, Newton trajectories, gentlest ascent dynamics, and that of the nudged elastic band or string variants) are all examples of paths that connect a minimum to a saddle. Sometimes, they coincide with reaction paths: continuous curves that monotonically ascend towards a saddle and then monotonically decrease to the next minimum.~\cite{bofill2020} (While we use the definition of the reaction path as described in Ref.~\onlinecite{bofill2020}, we note that there is disagreement in the literature surrounding the formal definition of a reaction path. For example, Ref.~\onlinecite{Heidrich1995} uses the term reaction path as a synonym for a curve in which all points of the reaction path sit in a valley. Further complicating matters, the path described by Ref.~\onlinecite{Heidrich1995} is sometimes referred to as the minimum energy path. We do not use this definition for a reaction path nor a minimum energy path in this work.) Curves that display turning points in $L$, and thus switch between increasing and decreasing in energy (i.e. curves that ``meander'') are not reaction paths. Here, we briefly compare the paths, which differ in  mathematical formulation and how they ascend to the saddle. Our comparison is done in Euclidean space.

The path of the nudged elastic band (NEB)~\cite{Jonsson1998} and string variants,~\cite{Weinan2002, behn2011, chaffey2012} converge to the steepest descent path from each side of a saddle point.~\cite{sheppard2011} This curve is often—but not always—referred to in literature as the minimum energy path (MEP) as it marks the ``path of least resistance'' or most probable reaction path between a chemical reactant and product.~\cite{dunitz1975} (An alternative definition for the MEP is a reaction path in which all points sit in a valley. This alternative definition has been well studied by Quapp and Bofill in Refs.~\onlinecite{ bofill2020, hirschquapp2004, bofill2012, bofill2015}. We do not use this definition in this work.) The tangent of this curve is parallel to the gradient except at fixed points. The discovery of the MEP using NEB requires an estimated path initialization and knowledge of two minimum,~\cite{palenik2021,halgren1977} which prevent it from being used with our algorithmic framework. Nevertheless, we present it here for comparison due to its significance.

A Newton trajectory (NT), also known as an isocline or reduced gradient following curve, is a curve in which the normalized vector field at each point along the curve is equal to a chosen, fixed vector.~\cite{quapp1998,hirsch2004} Unlike the path of the NEB, it carries no physical interpretation. A gradient extremal, as previously described, is a curve comprised of points in which the gradient is always an eigenvector of the Hessian of the potential. A gradient extremal does, however, bear physical significance as it usually follows a valley floor or ridge. Gentlest ascent dynamics forms a new dynamical system that is designed to search for and follow the eigenvector corresponding to the smallest eigenvalue of the Hessian of the potential such that a trajectory ascends towards a saddle point in the least steep manner possible.~\cite{Weinan2010} Note that this trajectory is not guaranteed to converge to a saddle. A GAD curve is different than a gradient extremal, due to the presence of the ascent term unless the gradient extremal is confluent to an MEP. Further comparisons of gradient extremals, GAD, and NTs can be found in Refs.~\onlinecite{hirsch2004, hirschquapp2004, quapp2003,bofill2013}.

\begin{figure}
\centering
\begin{subfigure}[]
  \centering
  \includegraphics[width=1\linewidth]{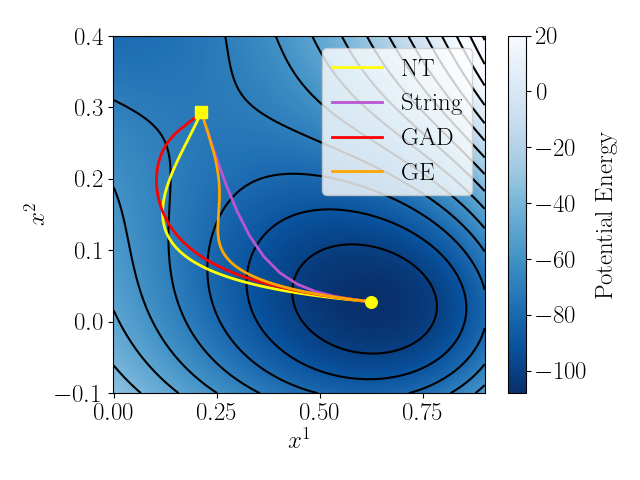}
  \label{fig:euc_comp_mb}
\end{subfigure}%
\begin{subfigure}[]
  \centering
  \includegraphics[width=1\linewidth]{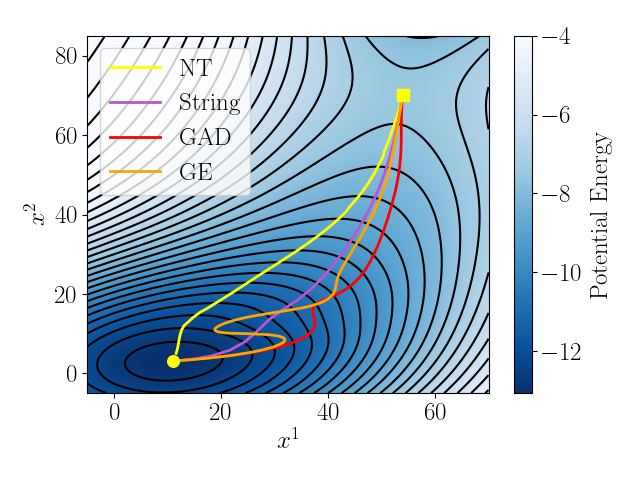}
  \label{fig:euc_comp_mb1}
\end{subfigure}
\caption{\textbf{(a)} Gradient extremal, the path of NEB (calculated via the String method), GAD (initialized with a velocity equivalent to the smallest eigenvector of the Hessian of the potential), and NT (with a horizontal slope) between the rightmost minimum and nearby saddle on the M\"{u}ller-Brown potential, which is further defined in~\eqref{eq:muller-brown}. \textbf{(b)} The four curves starting at the bottom minimum of the potential defined by~\eqref{eq:yannik}. Circles mark minima; squares mark saddles. All four curves are unique.}
\label{fig:Euc_Comparison}
\end{figure}

Figure~\ref{fig:Euc_Comparison}
illustrates the four curves on the M\"{u}ller-Brown potential and~\eqref{eq:yannik}. We emphasize that the different techniques do produce four distinct curves, and that the presence of ``meandering'' depends on the potential energy surface.\footnote{See~\cite{bofill2020} for a mathematical description of when NT, GAD, and GE exhibit turning points.} The gradient extremal curve is the most direct route (aside from that of the string method) to the saddle in the M\"{u}ller-Brown potential, but is the most indirect on~\eqref{eq:yannik} due to the presence of turning points. We have previously demonstrated the use of Newton trajectories and GAD within the algorithmic framework presented here to locate saddle points on Riemannian manifolds adaptively revealed in Refs.~\onlinecite{Bello-Rivas2023,bellorivas2023gentlest} respectively. The selection of which path is optimal for a given problem hinges on the shape of the potential energy surface, which, in the context of the algorithm, is likely unknown \textit{a priori}. 

\section{Numerical Example}
We demonstrate the algorithm on the 2D M\"{u}ller-Brown potential surface
\begin{align}
    \label{eq:muller-brown}
U(w^1, w^2)=\sum_{i = 1}^4 A_i \, \exp&\left( a_i (w^1 - w^1_{0i})^2 \right. \nonumber \\ 
&\left. + b_i (w^1 - w^1_{0i}) (w^2 - w^2_{0i}) \right. \nonumber \\
&\left. + c_i (w^2 - w^2_{0i})^2\right)
\end{align}

mapped onto the unit sphere using the transformation
\begin{equation*}
V(x,y,z) = (U \circ \kappa)(\arctan{(y/x)}, \arctan{(z/\sqrt{x^2+y^2}})
\end{equation*}
where $\kappa(k^1, k^2)=(1.973521294k^1-1.85, 1.750704373k^2+0.875)$ is an affine mapping. The constants in~\eqref{eq:muller-brown} are listed in Table~\ref{tab:mb}.

\begin{figure}
\centering
\begin{subfigure}[]
  \centering
  \includegraphics[width=0.4\linewidth]{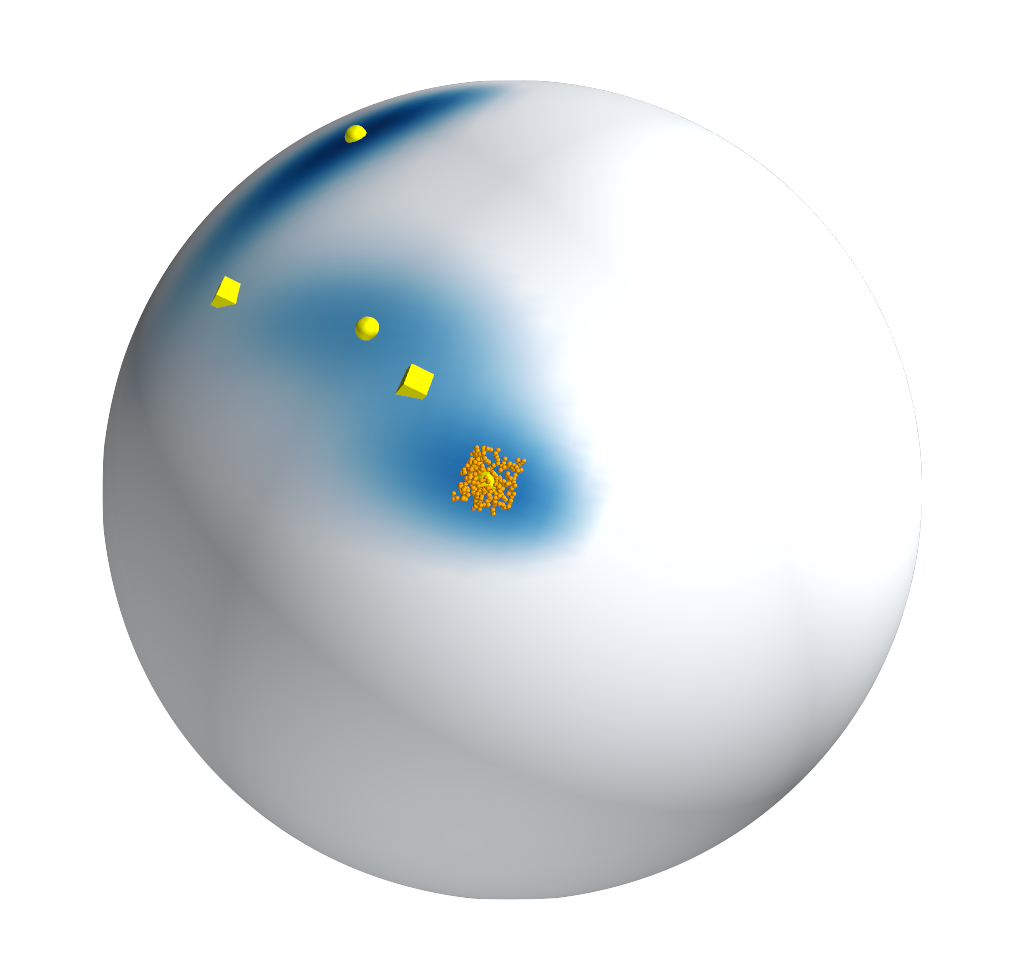}
\end{subfigure}%
\begin{subfigure}[]
  \centering
  \includegraphics[trim={13.5cm 14cm 1cm 1cm},clip,width=0.4\linewidth]{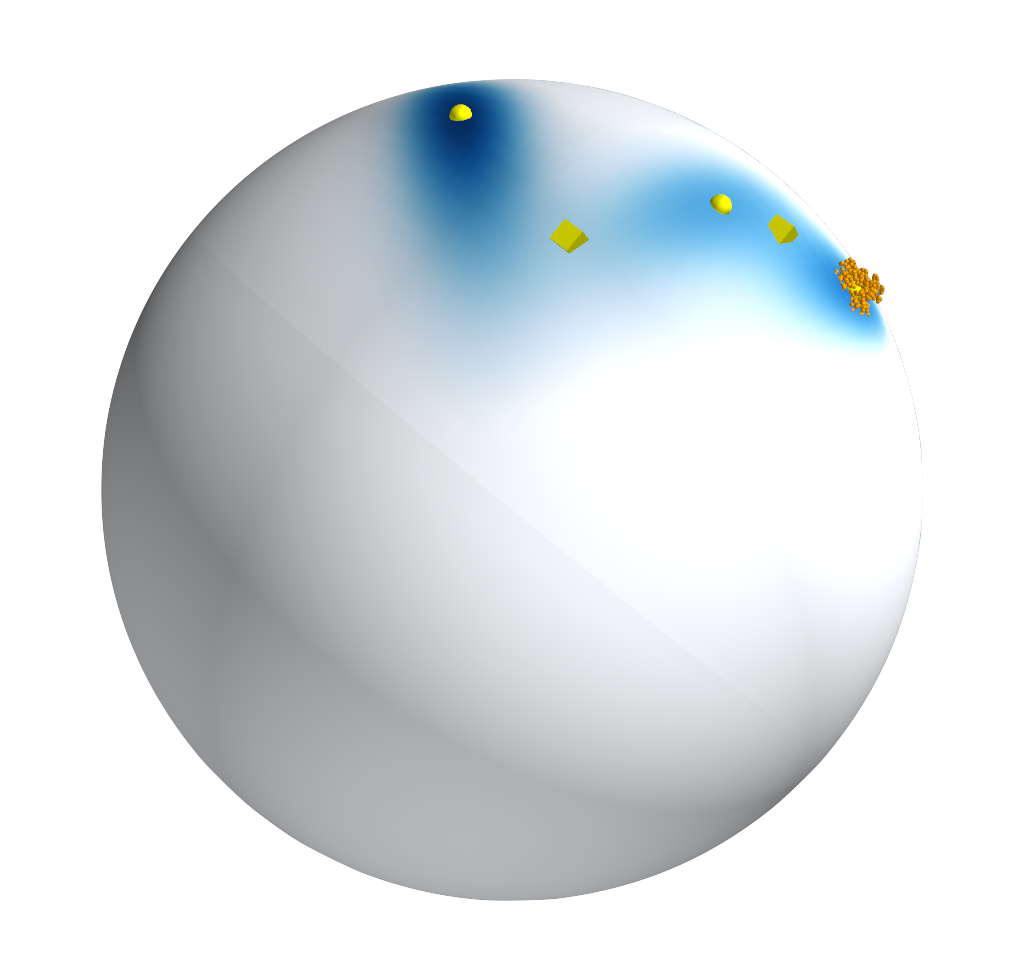}
\end{subfigure}
\begin{subfigure}[]
  \centering
  \includegraphics[width=0.50\linewidth]{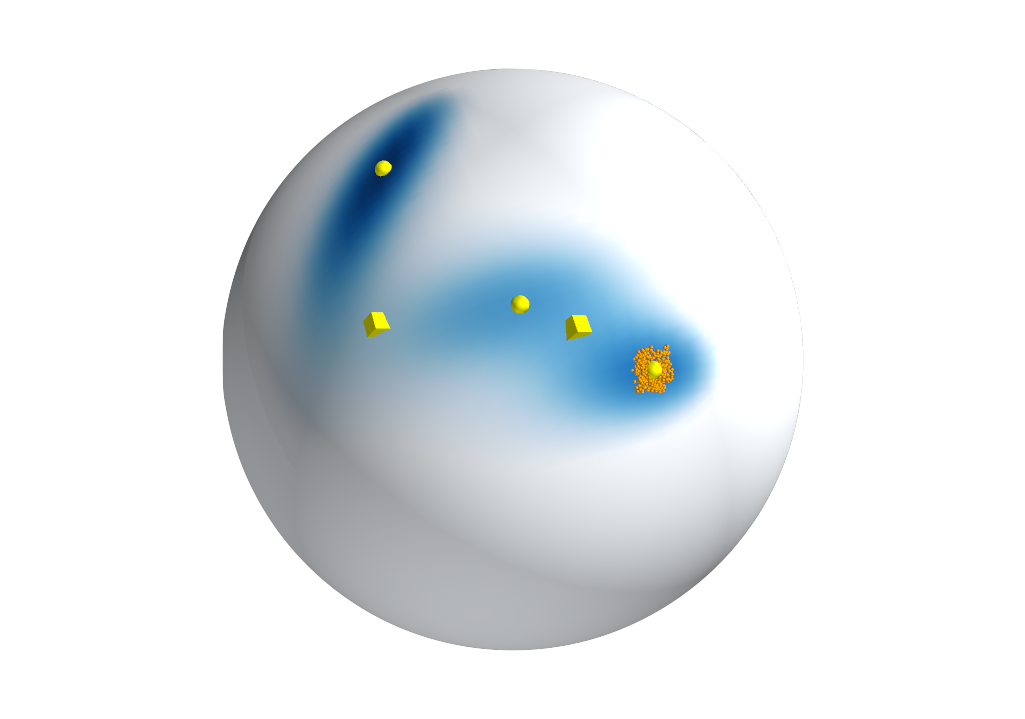}
\end{subfigure}%
\begin{subfigure}[]
  \centering
  \includegraphics[trim={13.5cm 10cm 1cm 2.5cm},clip,width=0.47\linewidth]{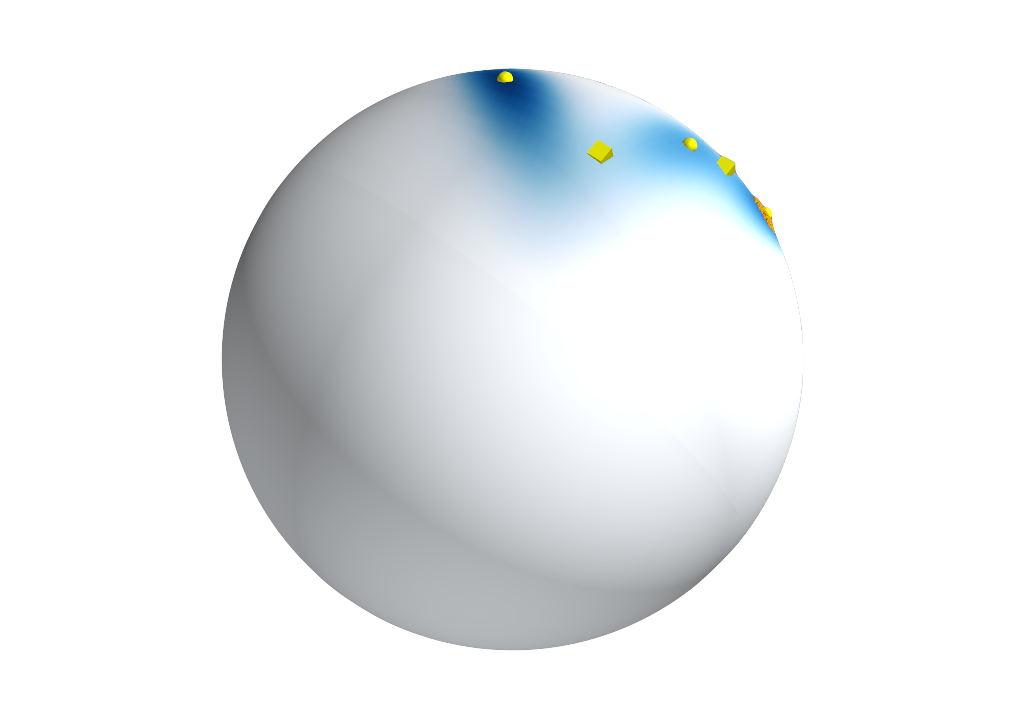}
\end{subfigure}
\caption{500 samples drawn around the rightmost minimum of the M\"{u}ller-Brown potential mapped onto a unit sphere. \textbf{(a,b)} Samples are first drawn by randomly perturbing the known minima in ambient space. \textbf{(c,d)} The trajectories are then run for a short time such that they have fallen onto the manifold. These final samples, now on the manifold, are used in the subsequent steps.}
\label{fig:Sampling}
\end{figure}

\begin{figure}
\centering
\begin{subfigure}[]
  \centering
  \includegraphics[width=0.9\linewidth]{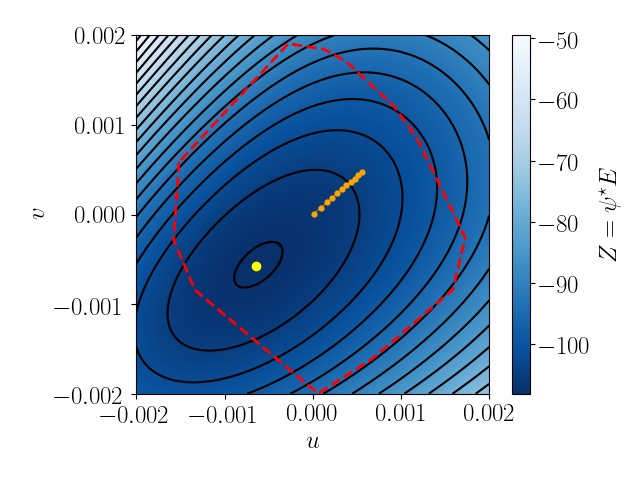}
  \label{fig:sub1}
\end{subfigure}%
\begin{subfigure}[]
  \centering
  \includegraphics[trim={4cm 4cm 4cm 4cm},clip,width=0.9\linewidth]{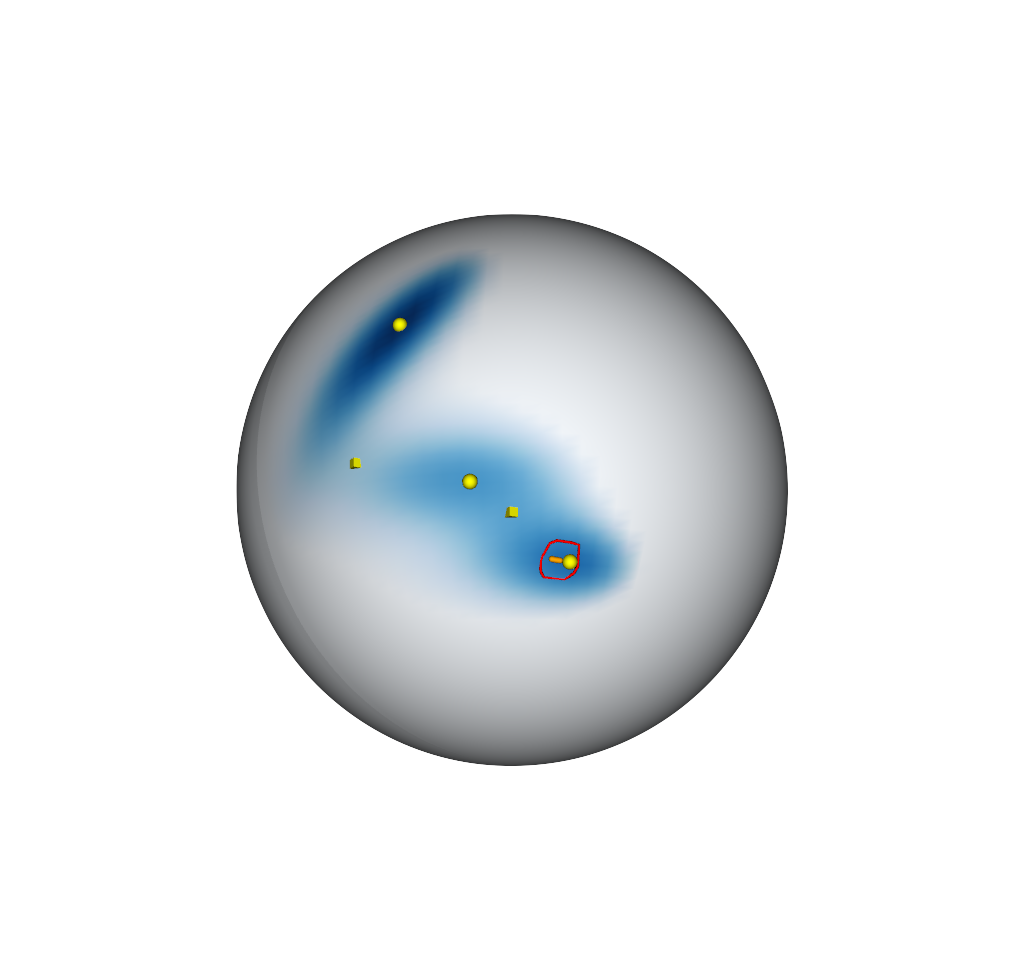}
\end{subfigure}
\caption{The second iteration of the algorithm beginning from the rightmost minimum of the M\"{u}ller-Brown potential mapped onto a unit sphere. \textbf{(a)} depicts the gradient extremal resolved in the learned latent chart. \textbf{(b)} shows the same path mapped back onto the sphere.  The red outline depicts the boundary of our learned latent chart.}
\label{fig:MB_step1}
\end{figure}

\begin{figure}
\centering
\begin{subfigure}
  \centering
  \includegraphics[width=0.9\linewidth]{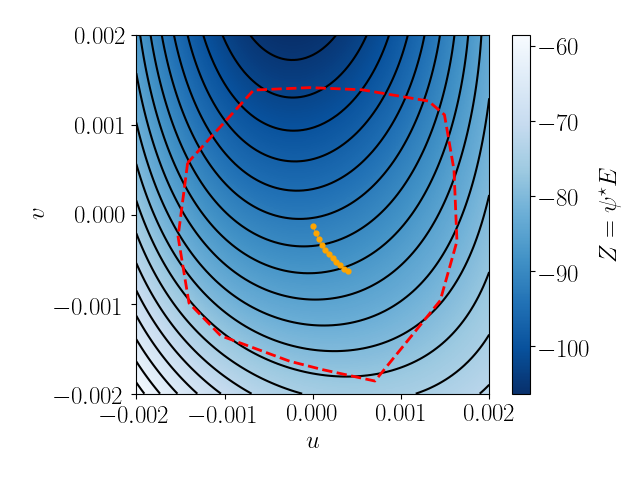}
  \label{fig:subc}
\end{subfigure}%
\begin{subfigure}
  \centering
  \includegraphics[trim={4cm 4cm 4cm 4cm},clip,width=0.9\linewidth]{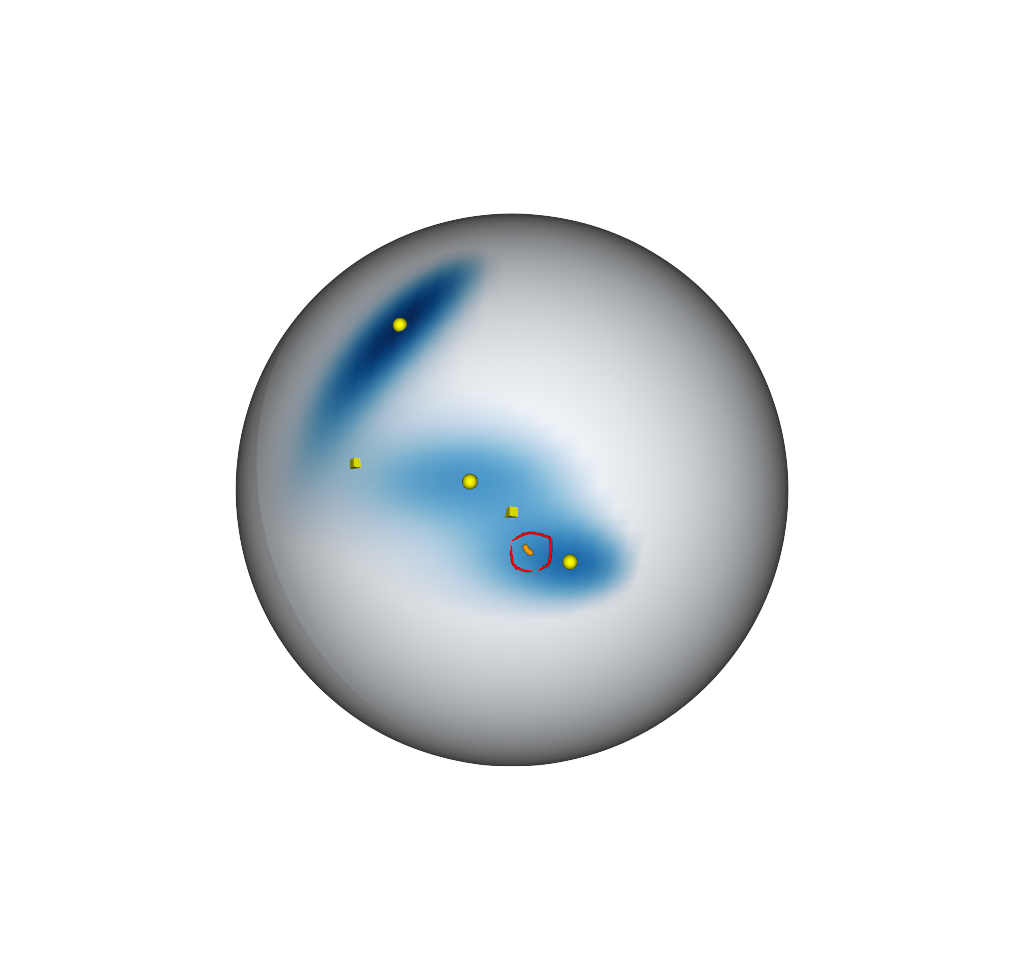}
\end{subfigure}
\caption{The fifth iteration of the algorithm beginning from the rightmost minimum of the M\"{u}ller-Brown potential mapped onto a unit sphere. \textbf{(a)} depicts the gradient extremal resolved in the learned latent chart. \textbf{(b)} shows the same path mapped back onto the sphere.  The red outline depicts the boundary of our learned latent chart.}
\label{fig:MB_step5}
\end{figure}

\begin{figure}
\centering
\begin{subfigure}
  \centering
  \includegraphics[trim={18cm 10cm 18cm 12cm},clip,width=0.75\linewidth]{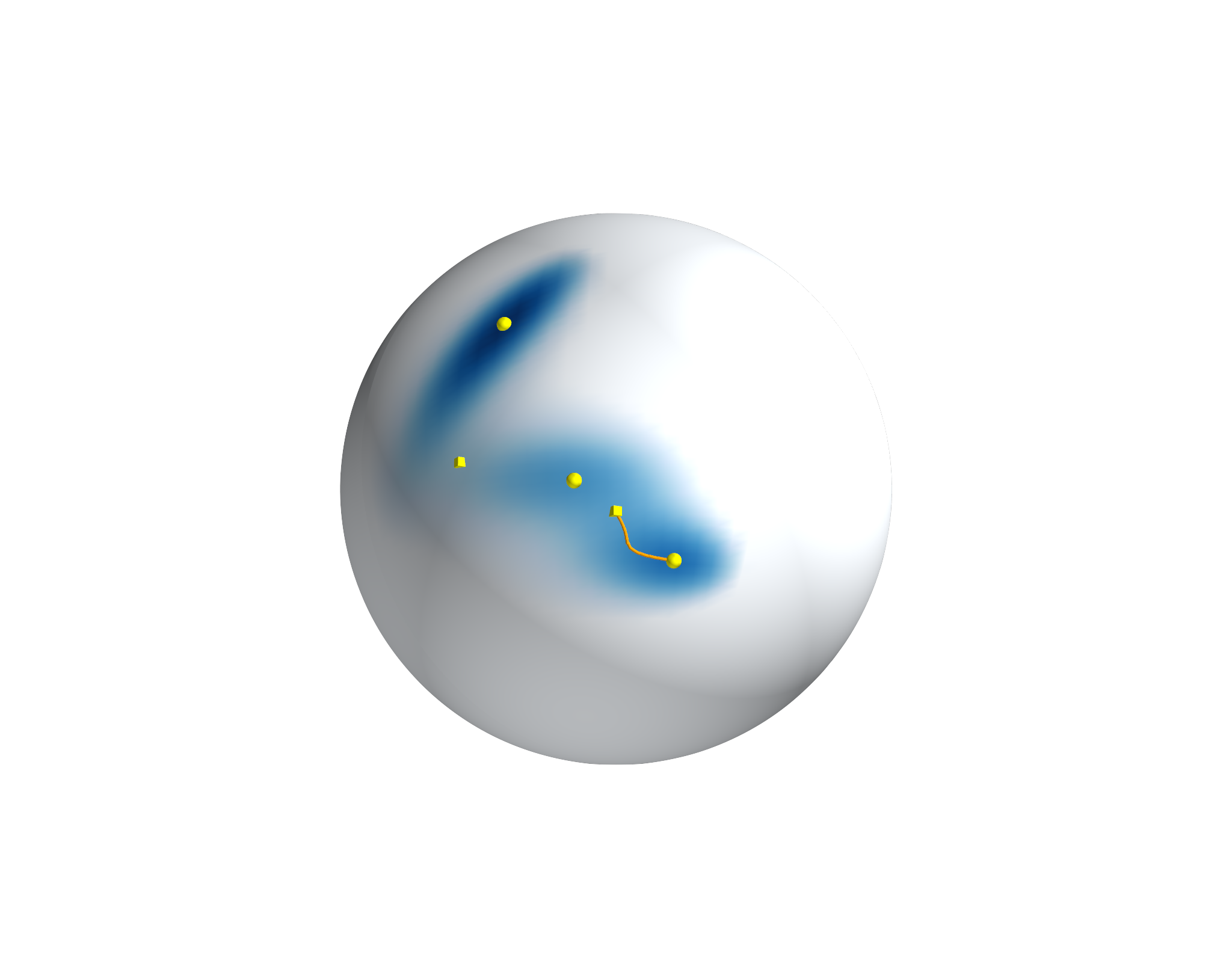}
\end{subfigure}%
\begin{subfigure}
  \centering
  \includegraphics[trim={8cm 6cm 8cm 6cm},clip,width=0.75\linewidth]{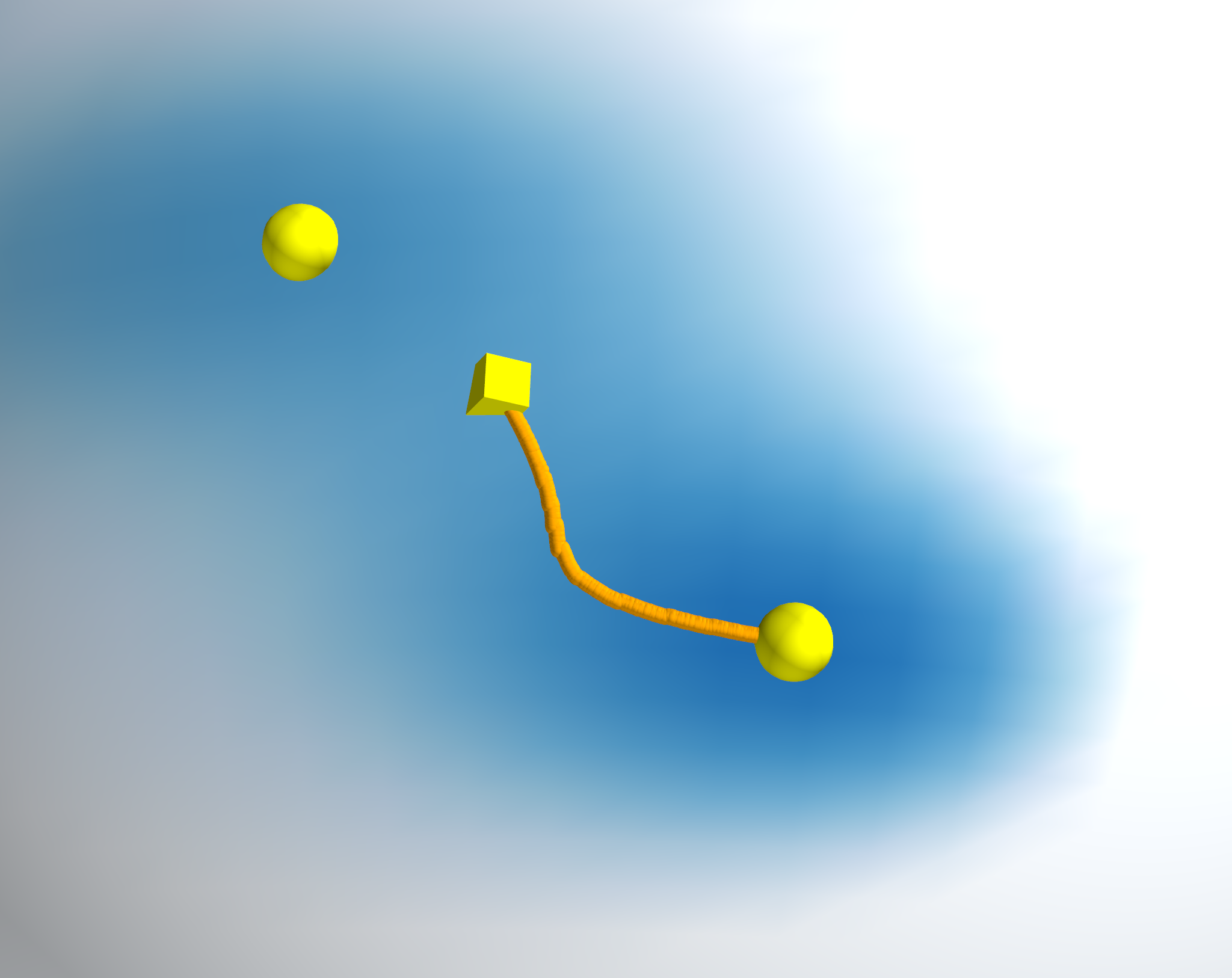}
\end{subfigure}
\caption{The resolved gradient extremal after 11 iterations successfully reaches the saddle point.}
\label{fig:MB_final}
\end{figure}

\begin{table}
\caption{\label{tab:mb}Coefficients of the planar M\"{u}ller-Brown potential.}
\begin{ruledtabular}
\renewcommand{\arraystretch}{1.25}
\begin{tabular}{c c c c c c c}
 $i$ & $A_i$ & $a_i$ & $b_i$ & $c_i$ & $y^1_{0i}$ & $y^2_{0i}$\\ [0.5ex]
 \hline
 1 & $-200$ & $-1$ & 0 & $-10$ & 1 & 0\\
 2 & $-100$ & $-1$ & 0 & $-10$ & 0 & 0.5 \\
 3 & $-170$ & $-6.5$ & 11 & $-6.5$ & $-0.5$ & $-1.5$ \\
 4 & 15.0 & 0.7 & 0.6 & 0.7 & $-1$ & 1\\ [1ex]
\end{tabular}
\end{ruledtabular}
\end{table}
Beginning at the rightmost minimum, we sample 500 points in the neighborhood of an initial point as in Figure \ref{fig:Sampling}, 
learn a local chart using diffusion maps, and learn the mappings between ambient space and charts using Gaussian processes as explained in Section~\ref{sec:Algorithm}. Figure~\ref{fig:MB_step1}
shows the second iteration of a gradient extremal in the learned latent space and its equivalence on the sphere in ambient space; Figure~\ref{fig:MB_step5} 
shows the fifth. We reached the saddle point after iterating through a sequence of 11 successive charts; the full gradient extremal can be seen in Figure~\ref{fig:MB_final}.
This example and a generic implementation of our algorithm can be found at \href{https://github.com/tasiag/gradient-extremals-on-manifolds}{https://github.com/tasiag/gradient-extremals-on-manifolds}.

\section{Conclusion}
We have presented an algorithmic framework to locate saddle points (and other fixed points) of vector fields on manifolds gradually revealed by point clouds. The methodology assumes that a low-dimensional manifold, which can be uncovered on-the-fly by manifold learning techniques, quickly attracts trajectories of the high-dimensional dynamical system. Gradient extremals, curves that are known to connect fixed points, can then be followed on the manifold to the next point of interest. The technique does not require \textit{a priori} knowledge of the manifold, nor of good collective variables; rather, it requires only an initial starting point close to a minimum and the ability to sample the dynamical system. The technique does not exhaustively search the whole space, but predicts a single new location to sample at, so that the simulator is efficiently and effectively biased towards the next point of interest.

We demonstrate the technique on the M\"{u}ller-Brown potential mapped to a sphere, without prior knowledge of the manifold. The technique is applicable both to (effective) gradient problems, where it helps discover transition states, as well as to general vector fields, where it helps connect isolated solution branches. We also compare how a gradient extremal curve compares to our prior work, which shares algorithmic elements but uses Newton trajectories~\cite{Bello-Rivas2023} and GAD curves~\cite{bellorivas2023gentlest} to locate saddle points. We expect the collection of the three techniques to act as a robust toolkit that can scale to systems of higher complexity.

\begin{acknowledgments}
We thank Yannick De Decker (Universit\'{e} Libre de Bruxelles, Belgium) for providing \eqref{eq:yannik}.

The work of A.S.G, J.M.B.-R., and I.G.K.  has been partially supported by the US Air Force Office of Scientific Research (AFOSR MURI) and the US Department of Energy. The work of A.S.G. has been partially supported by T32CA153952 from NCI. The work of H.V. has been supported by grant 1179820N from the FWO (Belgium).
\end{acknowledgments}

\section*{Data Availability Statement}

The data that support the findings of this study are openly available at the repository gradient-extremals-on-manifolds on Github,  \href{https://zenodo.org/record/8384160}{DOI 10.5281/zenodo.8384160}.

\nocite{*}
\bibliography{aipsamp}

\end{document}